\documentclass[letterpaper, 10 pt, conference]{ieeeconf}  

\IEEEoverridecommandlockouts                  
\usepackage{xspace,amssymb,epsfig,syntonly}
\usepackage{epsfig,amsmath,color}
\usepackage{xspace,syntonly,empheq}
\usepackage{url}
\usepackage{mathtools}
\usepackage{cite}
\usepackage{graphicx}
\usepackage{algorithm,algorithmic}
\usepackage{verbatim}
\usepackage{amssymb}
\usepackage{amsthm}
\usepackage{enumerate}
\usepackage{nicefrac}

\usepackage{subfig}

\renewcommand{\hat}{\widehat}
\renewcommand{\tilde}{\widetilde}

\newcommand{\cP}{{\cal P}}

\newcommand{\cX}{{\cal X}}

\newcommand{\reals}{\mathbb{R}}

\renewcommand{\qed}{\hfill\blacksquare}

\usepackage{epstopdf}

\newtheorem{definition}{Definition}
\newtheorem{proposition}{Proposition}

\newtheorem{theorem}{Theorem}
\newtheorem{corollary}{Corollary}
\theoremstyle{definition}
\theoremstyle{definition}
\newtheorem{remark}{Remark}
\newtheorem{observation}{Observation}

\title{\LARGE \bf
Inner Approximation of Minkowski Sums: A Union-Based Approach and Applications to Aggregated Energy Resources
}

\author{Md Salman Nazir, Ian A. Hiskens, Andrey Bernstein, and Emiliano Dall'Anese 
\thanks{This work was supported  in part by the National Science Foundation through grant ECCS-1810144, and by the Laboratory Directed Research and Development Program at the National Renewable Energy Laboratory (NREL).}
\thanks{Md Salman Nazir and Ian A. Hiskens are with the department of Electrical Engineering and Computer Science, University of Michigan, Ann Arbor, MI. Andrey Bernstein and is with NREL, Golden, CO. Emiliano Dall'Anese is with the department of Electrical, Computer, and Energy Engineering, University of Colorado Boulder, Boulder, CO.}
}

\begin{document}

\maketitle
\thispagestyle{empty}
\pagestyle{empty}

\begin{abstract}
This paper develops and compares algorithms to compute inner approximations of the Minkowski sum of convex polytopes. As an application, the paper considers the computation of the  feasibility set of aggregations of distributed energy resources (DERs), such as solar photovoltaic inverters, controllable loads, and storage devices. To fully account for the heterogeneity in the DERs while ensuring an acceptable approximation accuracy, the paper leverages  a union-based computation and advocates homothet-based polytope decompositions. However, union-based approaches can in general lead to high-dimensionality concerns; to alleviate this issue, this paper shows how to define candidate sets to reduce the computational complexity. Accuracy and trade-offs are analyzed through numerical simulations for illustrative examples. 

\end{abstract}


\section{Introduction}\label{sec:Intro}

Power systems are in the process of accommodating an increased amount of distributed energy resources (DERs) -- solar photovoltaic (PV) systems, energy storage systems, and controllable demand-side resources just to mention a few. The flexibility from DERs can be leveraged to alleviate a number of operational challenges in the power grid \cite{Callaway2011, Hao2013, Muller2015} -- for example, to address  voltage regulation issues -- and to aid system-level operations by realizing the emerging vision of virtual power plants. To address this, key is to characterize  the aggregate flexibility from DERs~\cite{Bernstein2016,Kundu2018,Barot2017,Zhao2017,Muller2015,Muller2017,Alizadeh2014, Nazir2017, Hre2017}.

A general framework for characterizing DER flexibility is presented in \cite{Hao2015}, where methods to compute the aggregate flexibility using the Minkowski sum (M-sum) are also described.  M-sum can be computed accurately by summing all the vertices of given polytopes \cite{BrunnMink, Lien2007}; however, such approaches are not computationally feasible due to the exponential growth in complexity for large number of devices \cite{Barot2017, Muller2015, Weibel2007}. Hence, several works in the literature sought efficient algorithms to compute the M-sum. In \cite{Barot2017}, the authors provide an algorithm to compute an outer approximation of the M-sum; however, outer approximations might include infeasible points, which is undesirable especially if utilized in optimization settings. The authors in \cite{Muller2015} present an algorithm to compute an inner approximation of the M-sum by using zonotopes. Zonotopes have also been used for computing M-sum widely in literature due to their features which allow easily summing them to obtain the M-sum \cite{Althoff2011, Fukuda2004, Muller2015, Muller2017}. One major limitation is that zonotopes are centrally symmetric objects, hence if original resource polytopes are not symmetric (as will be shown for the inverter case), approximating those using single symmetric polytopes might lead to a conservative estimate of the aggregate flexibility. 
 
The authors in \cite{Zhao2017} present an algorithm to compute both inner and outer approximations of M-sum by using homothets. Given a prototype set, it can be scaled and translated to fit inside (for inner approximation) or just outside (for outer approximation) of given resource polytopes. However, choosing an arbitrary prototype shape may lead to a conservative estimate of the aggregate flexibility. 


This paper aims to extend zonotope- and homothet-based approaches for computing M-sums. The focus is on finding inner approximations of the M-sum so that the feasibility of control actions is guaranteed. The flexibility provided by inverter-interfaced devices as well as controllable loads is considered. 
To that end, a polytopic representation of the feasible operating region of an inverter is presented first. For certain special cases, we  provide  analytical expressions for the inner approximation of the M-sum by leveraging homothet-based representations.
However, as the level of heterogeneity increases, these analytical expressions might provide highly conservative estimates. Therefore, we propose to employ multiple homothets (here, axis-aligned boxes, which are essentially zonotopes \cite{Muller2015}) per device polytope, and show how to efficiently perform the M-sum computation. 

Our approach consists of: (i) a decomposition procedure to find a number of homothet-based sub-polytopes per device polytope; and (ii) performing the M-sum computation from the union of such sub-polytopes. We provide  asymptotic guarantees on the accuracy of the approximation, which is generally difficult to achieve for non-vertex based M-sum algorithms. Since the number of unions grows exponentially \cite{Tiwary2008, Weibel2007} with the number of devices and the number of sub-polytopes per device, techniques to limit the computational complexity of the methods are explored. The accuracy versus complexity trade-offs are investigated. 

The rest of the paper is organized as follows. Section \ref{sec:InvFlexibility} provides an overview of flexible operating region of inverter-interfaced devices and controllable loads, such as pumps, variable speed drives, electric vehicles (EVs) and thermostatically controlled loads (TCLs). A discretization technique to obtain a convex flexibility polytope is also presented. Section \ref{sec:AggrFlexM-sum} describes homothet-based approach to obtain the M-sum and proposes simple analytical expressions for its inner approximation. Section \ref{sec:UnionM-sum} presents a union-based M-sum algorithm, along with homothet-based decomposition technique. Section \ref{sec:Results} illustrates the effectiveness and accuracy of our techniques through numerical results. Finally, Section \ref{sec:Conclusion} concludes. 

\section{Flexibility Characterization}\label{sec:InvFlexibility}

We next focus on the characterization of the flexibility regions of inverter-based devices and controllable loads. 

\subsection{Inverter Feasible Set}\label{sec:InvSet}

\begin{figure}[t!]
	\begin{center}
		\includegraphics[scale=0.25]{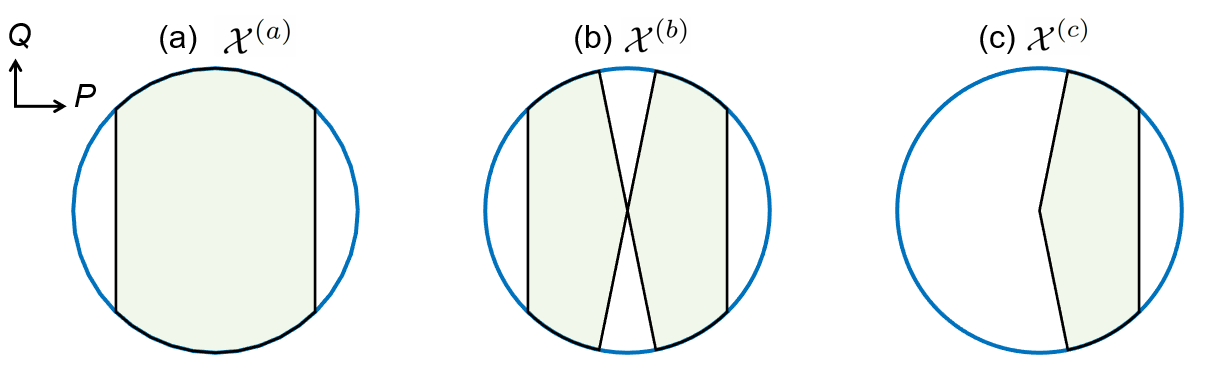} 
		\caption{Feasible sets of inverters. (a) $\cX(S,\underline{P},\overline{P})$, (b) $\cX(S,\underline{P},\overline{P},\theta_{})$, and (c) $\cX(S,\overline{P},\theta)$. }
		\label{fig:InverterSets}
	\end{center}	
\end{figure}

Let $\cX \subseteq \reals^2$ be the set that contains the inverter's real and reactive power operating points, $x = [P_{}, Q_{}]^T, x \in \reals^2$. Then, $\cX$ can be written as (see, e.g., \cite{opfPursuit}) 
\begin{align} \label{eq:inv1}
\cX(S,\underline{P},\overline{P}) = \{(P,Q): S\underline{P} \le  P  \le  S\overline{P},  Q^2 \le S^{2} - P^2 \}.
\end{align}
Here, for a PV system, $S$ is the apparent power rating of the inverter, $\underline{P} = 0$ and $\overline{P}$ (normalized w.r.t $S$) is the available power based on solar irradiance. For a storage device interfaced with an inverter, $S$ is the inverter's rating; $\underline{P}$ and $\overline{P}$ (normalized w.r.t $S$) are the minimum and maximum real power available at a specific time. Note that $\cX(S,\underline{P},\overline{P})$ is a convex set. Additionally, to enforce a minimum power factor of $\text{cos}(\theta_{})$, the following constraint can be included:
\begin{align}\label{eq:inv3}
| Q_{} | \le \text{tan}(\theta_{}) P_{};  
\end{align}
let $\cX(S,\underline{P},\overline{P},\theta_{})$ denote the resulting set, and note that it is generally non-convex. However, for PV systems, since $\underline{P} = 0$, it is convex; with a slight abuse of notation we let $\cX(S,\overline{P},\theta)$ denote this latter set.
Fig. \ref{fig:InverterSets} illustrates these feasible sets - $\cX(S,\underline{P},\overline{P})$ in (a), 
$\cX(S,\underline{P},\overline{P},\theta_{})$ in (b) and $\cX(S,\overline{P},\theta)$ in (c). For the rest of this paper, we will only focus on the convex cases (a) and (c). 

\subsection{Inverter Flexibility Polytope}\label{sec:InvPolytope}

Polytopes can generally be expressed with vertices (V-rep) or with half-space constraints (H-rep). The H-rep is useful for optimization purposes \cite{Barot2017, Zhao2017}. Hence, a polytopic representation of $\cX$ will be developed in this section.  
\begin{definition}
Let $\cP_{} = \{ A_{} x \le b_{} \}$, where $ x  = [P_{}, Q_{}]^T \in \reals^2$, $A_{} \in \reals^{m \times 2}$, and $b \in \reals^{m}$. If $\cP_{} \subseteq {\cX}$, then $\cP_{}$ is an \emph{inner approximation} to ${\cX_{}}$. 
\end{definition}
To obtain $A_{}$, $b_{}$, first we inscribe an $N$-sided polygon inside the circle $P_{}^2 + Q^{2} = S_{}^2$. Assume $N$ is even and $N \ge 4$. The angle formed between two successive vertices of the polygon, $\alpha$, can be found as, $\alpha = \frac{2 \pi}{N}$

The set of vertices $\{(P_{j}, Q_{j})\}_{j=1}^N$ can be found as
\begin{align}
P_{j} &= S \text{cos}\big((j-1)\alpha \big), \quad j = 1,2,...,N, \\
Q_{j} &= S \text{sin}\big((j-1)\alpha \big), \quad j = 1,2,...,N.
\end{align}
Additionally, by convention, $P_{N+1} = P_1$ and $Q_{N+1} = Q_1$.
From these vertices, the slopes $m_{j}$ can be computed as,
\begin{equation}\label{eq:mslope}
m_{j}  = \frac{(Q_{j+1} - Q_{j})}{(P_{j+1} - P_{j})}  , \quad j = 1,2,...,N.
\end{equation}
Then, the constraint set for the H-rep of $Q_{}^2 \le S^{2}_{} - P_{}^2$ can be obtained as, 
\begin{align}
 (Q_{} - Q_{j}) & \le m_{j} (P_{} - P_{j}) , \quad j = 1,2,...,\frac{N}{2},  \label{eq:circlconstr1}  \\
 -(Q_{} - Q_{j}) & \le -m_{j} (P_{} - P_{j}) , \quad j = \frac{N}{2}+1,...,N.  \label{eq:circlconstr2}
\end{align}

Overall, for $\cX(S,\underline{P},\overline{P})$, we obtain the polytope 
\begin{align} \nonumber  
\cP(S,\underline{P},\overline{P}) = \{(P,Q): S\underline{P} \le  P_{}  \le  S\overline{P}, \, \eqref{eq:circlconstr1} \text{ and } \eqref{eq:circlconstr2} \}   
\end{align}
Similarly, for $\cX(S,\overline{P},\theta)$, we obtain
\begin{align} \nonumber  
\cP(S,\overline{P},\theta) = \{(P,Q): 0 \le  P \le  S\overline{P}, \, \eqref{eq:inv3}, \, \eqref{eq:circlconstr1} \text{ and } \eqref{eq:circlconstr2} \}   
\end{align}

\begin{figure}[t!]
	\begin{center}
		\includegraphics[scale = 0.27]{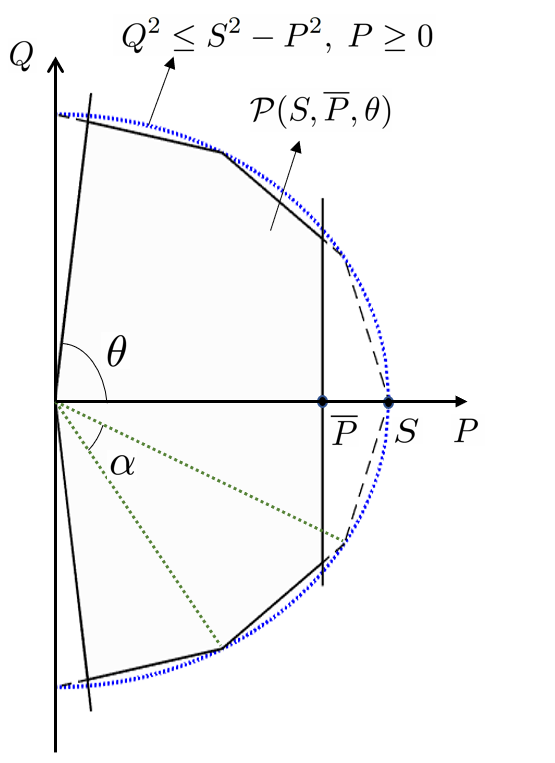} 
		\caption{Flexibility polytope for a photovoltaic inverter}
		\label{fig:InverterPolytope}
	\end{center}	
\vspace{-.4cm}
\end{figure}

The H-reps of $Q_{}^2 \le S^{2}_{} - P_{}^2$ and $\cP(S,\overline{P},\theta)$ are illustrated in Fig. \ref{fig:InverterPolytope}. In this case, since $\cP(S,\overline{P},\theta)$ spans two-quadrants, only the positive half-circle is linearized to avoid redundancy.

Note that, by construction, $\cP_{}$ is a convex polytope and for any finite $N$, $\cP_{} \subset \cX_{}$. Let $A_{\cP}$ denote the area of $\cP$ and $A_{\cX}$ denote the area of the entire feasibility set $\cX$. Then, the area ratio, $\eta := \frac{A_{\cP}}{A_{\cX}}$.

\begin{proposition}\label{PropInv1}
Consider $\cX(S,\underline{P},\overline{P})$, with $\underline{P}=-1$, $\overline{P}= 1$, or $\cX(S,\overline{P},\theta)$ with $\overline{P}= 1$ and $\theta = \pi/2$.  Then, $\eta = \frac{sin\alpha}{\alpha}$ and as $N \rightarrow \infty$, $\eta \rightarrow 1$.
\end{proposition}

\begin{proposition}\label{PropInv2}
Consider $\cX(S,\underline{P},\overline{P})$, with $\underline{P}>-1$, $\overline{P}< 1$, or $\cX(S,\overline{P},\theta)$ with $\overline{P}< 1$ and $0 \le \theta < \pi/2$. Then, for $N$ sufficiently large, $\eta \approx 1$ .
\end{proposition}

For the proof of this and other results of this paper, see the Appendix.

Applying Proposition \ref{PropInv1} with $N = 6$ yields $\eta = 0.83$; with $N = 12$ yields $\eta = 0.95$; and with $N = 24$ yields $\eta = 0.99$. Hence, $N = 24$ will be used in our simulations below.

 \subsection{Flexibility from Controllable Loads}\label{sec:FlexLoads}

We also consider controllable loads, such as variable speed drives and pool pumps. Their feasible set ${\cX_{}}$, considering only real power, can simply be written as
\begin{align}
\cX(\underline{P},\overline{P}) = \{P: \quad  & \underline{P} \le  P_{L}  \le  \overline{P} \}. \label{eq:controllableload1} 
\end{align}

For storage-like DERs, such as TCLs and EVs, the internal states (e.g. temperature, state of charge) also evolve with time. Consider $M$ time-intervals, indexed by $k = 1,2,...,M$. Let $e_{k}$ be the DER's normalized energy state, i.e. its state of charge (SOC), and $P_k$ be the real power consumed by the DER at time $k$. Then, using a generalized battery model, the dynamics of $e_{k}$ can expressed as, $e_{k+1} = a e_{k} + \gamma P_{k}$, where $a \in (0,1]$ is the energy dissipation rate, and $\gamma$ is the charging efficiency \cite{Hao2015, Zhao2017}. Knowing the initial SOC, $e_{o}$, the overall feasible set $\cX$ can be written as \cite{Barot2017}, 
\begin{multline}
\cX(  \underline{P},\overline{P},a,\gamma,e_0) = \Big\{ [P_k]^{\textrm{T}}: \underline{P} \le  P_{k} \le  \overline{P},  \\
 0 \le a^k e_{o} + \sum_{t=1}^{k} a^{k-t} \gamma P_{k}  \le 1, k = 1,2,...,M \Big\}.  \label{eq:gbm3} 
\end{multline}
with $ x = [P_k]^{\textrm{T}} \in \reals^M$.
For more details on the polytopic representation  of \eqref{eq:gbm3}, readers can refer to \cite{Barot2017, Zhao2017}. 


\section{Aggregation by Minkowski Sum}\label{sec:AggrFlexM-sum}

In this section, we describe how the flexibility from DERs of the same type can be aggregated using M-sum. Consider a population of $n_d$ devices, with indices $i = 1, 2,...,n_d$. Let $\cX_i$ denote the feasible set of device $i$. The aggregate flexibility, ${\cX}^{\textrm{Aggr}}$, can be found by computing the M-sum of ${\cX_i}$s as,
\begin{equation}\label{eq:defnMsum}
{\cX}^{\textrm{Aggr}} := {\cX_1} \oplus {\cX_2} \oplus ... \oplus {\cX_{n_d}} = \oplus_{i=1}^{n_d} \cX_i. 
\end{equation} 
where $\oplus$ denotes the M-sum. For computation of M-sum, applying \eqref{eq:defnMsum} is, however, not efficient, especially when $n_d$ is large \cite{Weibel2007}. Hence, zonotopes \cite{Fukuda2004,Muller2015} and homothet-based \cite{BrunnMink,Zhao2017} approaches have been shown to be useful. The applicability of homothets for our case will be presented next.

\subsection{Homothets and Minkowski Sum}\label{sec:M-sumbyHomothets}
Given a compact convex set $\cX_{0}$, $\beta_i \cX_{0} + t_i := \{ x \in \reals^2: \, x = \beta_i \zeta + t_i, \, \zeta \in \cX_{0} \}$ is a homothet of $\cX_{0}$, where $\beta_i > 0$ is a scaling factor and $t_i$ is a translation factor \cite{BrunnMink}. $\cX_{0}$ can be referred to as a prototype set \cite{Zhao2017}. 

Homothets are useful for computation of M-sums due to the following property \cite{BrunnMink, Zhao2017},
\begin{equation}\label{eq:M-sumhomothet1}
 \oplus_{i=1}^{n_d} (\beta_i{\cX_{0}}+t_i) =  \sum_{i=1}^{n_d} \beta_i {\cX_{0}} + \sum_{i=1}^{n_d} t_i.
\end{equation}
Hence, if $(\beta_i{\cX_{0}}+t_i) \subseteq {\cX_i}$, then,
\begin{equation}\label{eq:M-sumhomothet2}
 \oplus_{i=1}^{n_d} (\beta_i{\cX_{0}}+t_i) \subseteq  \oplus_{i=1}^{n_d} \cX_i.
\end{equation}

For example, consider $\cX_i(\underline{P_i},\overline{P_i})$, defined in \eqref{eq:controllableload1}. 
Take $\cX_{0}=\cX(1,1)$ and find $\beta_i$s and $t_i$s for all $\cX_i$. Then, applying \eqref{eq:M-sumhomothet1},
the aggregate flexibility from $n_d$ controllable loads, is simply given by, $\cX^{\textrm{Aggr}}=\cX(\sum_{i=1}^{n_d}\underline{P_i},\sum_{i=1}^{n_d}\overline{P_i})$. 

\subsection{Special Cases: Aggregate Flexibility from Inverter-interfaced Devices}\label{sec:InvbyHomothets}

Under certain conditions, the properties \eqref{eq:M-sumhomothet1} and \eqref{eq:M-sumhomothet2} of homothets  lead to simple analytical expressions for the M-sum of inverter-interfaced devices.  

For example, assume only the rated power of the inverters vary, while $\underline{P_i}$ and $\overline{P_i}$ are homogeneous. This situation can appear commonly when a collection of inverters have different ratings, $S_i$, but undergo similar solar irradiance conditions, which could be due to their geographic proximity. Their aggregate flexibility can be obtained by Theorem \ref{THM-Svary}. 

\begin{theorem}\label{THM-Svary}
Consider $\cX_i(S_i,\underline{P_i},\overline{P_i})$, where $\underline{P_i} = \underline{P_0}$ and $\overline{P_i} = \overline{P_0}$ for $i = 1,...,n_d$. 
The aggregate flexibility set is then given by
$$\oplus_{i=1}^{n_d} \cX_i(S_i,\underline{P_i},\overline{P_i}) = \cX \left(\sum_{i=1}^{n_d} S_i, \underline{P_0}, \overline{P_0}\right).$$
\end{theorem}

\begin{corollary}\label{Corr-Svary}
Consider $\cX_i(S_i,\overline{P_i},\theta_i)$, where $\overline{P_i} = \overline{P_0}$ and $\theta_i = \theta_0$ for $i = 1,...,n_d$. The aggregate flexibility set is given by ${\cX}(\sum_{i=1}^{n_d} S_i, \overline{P_0}, \theta_0)$.
\end{corollary}

Next, consider heterogeneous $S_i$, $\overline{P_i}$ and $\theta_i$ for $\cX_i(S_i,\overline{P_i},\theta_i)$. In this case, Theorem \ref{THM-hetero} applies. 

\begin{theorem}\label{THM-hetero}
Consider $\cX_i(S_i,\overline{P_i},\theta_i)$, $i = 1,...,n_d$, where $S_i,\overline{P_i},\theta_i$ are heterogeneous. Let $S_0 = \underset{i}\min S_i$, $\overline{P_0} = \underset{i}\min \overline{P_i}$, and $\theta_0 = \underset{i}\min \theta_i$ ($0 \le \theta_i \le \pi/2$). Then, $$\cX \big(n_d S_0, \overline{P_0}, \theta_0 \big) \subseteq \oplus_{i=1}^{n_d} \cX_i(S_i,\overline{P_i},\theta_i).$$ 
Moreover, strict equality holds if and only if all the parameters are homogeneous. 
\end{theorem}



\textbf{{Example 1.}} Take $n_d = 100$ inverters characterized by $\cX(S_i,\overline{P_i},\theta_i)$. First, consider, $S_i = 1, \theta_i = 1.45$ rad, $\forall i$. $\overline{P}$ is distributed uniformly between 0.75 and 1. By Theorem 2, using $\min_i \overline{P_i}$, the prototype set and the approximate M-sum were obtained. The discretization procedure described in section \ref{sec:InvPolytope} was used to obtain polytopic representations for the $n_d$ inverters. Using MPT toolbox \cite{MPTtoolbox} the actual M-sum and its area were computed. The ratio of the area of approximate M-sum polytope to the area of the true M-sum polytope was $0.90$. Next, assume that $\overline{P_i}$ is distributed uniformly between $0.5$ to $1$. The area ratio was found to be $\approx 0.71$. Next, a population was considered with all parameters being heterogeneous, $S_i$ uniformly distributed within [0.75,1], $\overline{P_i}$ within [0.75,1] and $\theta_i$ within [1.27,$\pi$/2] rad. In this case, the area ratio drastically reduced to 0.29. Thus, as the level of heterogeneity increases, applying Theorem 1 and 2 may lead to very conservative estimates.$\qed$

The accuracy of the M-sum approximated by using homothets depends on the choice of the prototype set. Typically, as the level of heterogeneity increases, the accuracy worsens considerably. Hence, in this paper, we present an approach to consider multiple homothets  per $\cX_i$ and show that the inner approximation of the M-sum can approach the true M-sum. 

\section{Union-based Minkowski Sum}\label{sec:UnionM-sum}

For a heterogeneous population, the shapes of the flexibility sets may vary considerably. Hence, choosing a single prototype set $\cX_0$ may be limiting and result in a conservative estimate of the M-sum. To address this, in this section, we show how to decompose each polytope into a union of homothetic sub-polytopes. The M-sum can then be computed by applying the distributivity property of M-sum, as elaborated in the next sub-section.

Our union-based approach can also be motivated by the optimization applications as follows. 
Given a collection of convex compact subsets $\{\cX_{\omega}\}_{{\omega} = 1}^{n_{\Omega}}$ of $\reals^M$, consider the union $\cX := \bigcup_{{\omega} = 1}^{n_{\Omega}} \cX_{\omega}$. Let $f: \reals^M \rightarrow \reals$ be a convex function, and consider the optimization problem,
\begin{align}\label{eq:optBasic} 
\textrm{(P0)}\quad &\min_{x \in \cX}  \hspace{.2cm} f(x). 
\end{align}
It is clear that (P0) is equivalent to:
\begin{align}\label{eq:optUnion} 
\textrm{(P1)}\quad & \min_{ \omega \in \{1, \ldots, n_{\Omega}\}} \min_{x \in \cX_{\omega}} \hspace{.2cm} f(x). 
\end{align}
In this case, if $\cX$ represents the aggregate flexibility set, we avoid computing the overall M-sum; instead, we find the optimal solution from the candidate solutions obtained from solving multiple sub-problems. 

\subsection{Distributivity Property of Minkowski Sum}\label{sec:DistributivityMsum}

Let each $\cX_i$, the set for the $i$-th of the $n_d$ DERs, be expressed by $n_i$ sub-sets. Let $W_i = \Big\{ (i,j): j=1,...,n_i \Big\}$. Then,
\begin{align}
\cX_i &= \cup_{\omega \in W_i} \cX_\omega.
\end{align}
Also, let $\Omega$ be the Cartesian product of all $W_i$, i.e. $\Omega = W_1 \times ... \times W_{n_d}
= \Big\{(\omega_1,...,\omega_{n_d}) : \omega_i \in W_i, \forall i=1,...,n_d \Big\}$.
Then, by the distributivity property of M-sum \cite{BrunnMink},
\begin{align}
\oplus_{i=1}^{n_d} \mathcal{X}_i &= \oplus_{i=1}^{n_d} \Big(
\cup_{\omega \in W_i} \mathcal{X}_\omega \Big)  \label{eq:Distrib3} \\
&= \cup_{(\omega_1,...,\omega_{n_d})\in \Omega} \Big(
\oplus_{i=1}^{n_d} \mathcal{X}_{\omega_i} \Big). \label{eq:Distrib4}
\end{align}
Because for each of the $n_d$ DERs, one can choose from $n_i$ sub-sets, the cardinality of $\Omega$ is $n_{\Omega} = \prod_{i=1}^{n_d} n_{i}$. Note that while \eqref{eq:Distrib4} holds with equality, one may chose any number of subsets from $\Omega$ and obtain $\bar{\Omega} \subseteq \Omega$. Then, 
\begin{equation}\label{eq:Distrib5}
\cup_{(\omega_1,...,\omega_{n_d})\in \bar{\Omega}} \Big(
\oplus_{i=1}^{n_d} \mathcal{X}_{\omega_i} \Big)\subseteq  \oplus_{i=1}^{n_d} \cX_i,
\end{equation}
i.e. an inner approximation of $\oplus_{i=1}^{n_d} \cX_i$ is obtained.
For example, consider $\cX_1$ and $\cX_2$ and assume $n_1 = n_2 = 2$. Then, by \eqref{eq:Distrib4},  
\begin{multline}
\cX_1 \oplus \cX_2 = \big(\cX_{(1,1)}  \oplus \cX_{(2,1)}\big) \cup \big(\cX_{(1,1)} \oplus \cX_{(2,2)}\big) \\ \cup \big(\cX_{(1,2)} \oplus \cX_{(2,1)}\big) \cup \big(\cX_{(1,2)}  \oplus \cX_{(2,2)}\big).  
\end{multline}
Of course, $\big(\cX_{(1,1)} \oplus \cX_{(2,2)}\big) \subseteq \big(\cX_1 \oplus \cX_2\big)$.

Finally, expressing every subset of $\cX_i$  as a homothet of the same prototype set $\cX_0$, from \eqref{eq:Distrib4} and \eqref{eq:M-sumhomothet1} we obtain, 
 \begin{multline}\label{eq:MsumUnionHomothet}
 \cup_{(\omega_1,...,\omega_{n_d})\in \Omega} \Big(
\oplus_{i=1}^{n_d} \big(\beta_{\omega_i} \cX_0 + b_{\omega_i}\big) \Big)
 \subseteq  \oplus_{i=1}^{n_d} \cX_i.
 \end{multline}

The challenges associated with union-based M-sum include:
(1) Optimally partitioning a given polytope into convex sub-polytopes, 
(2) Analyzing the trade-offs between computational complexity and accuracy with increasing $n_{i}$ and $n_{d}$.
To efficiently handle these, a decomposition algorithm is proposed next. 


\subsection{Homothet-based Polytope Decomposition (HPD) }\label{sec:PDecomposition}

The key idea here is to decompose each of the given $M$-dimensional convex polytopes $\cP := \{x: \, Ax \leq b \}$ into a number of homothets. Consider axis-aligned boxes. Let the lower and upper boundaries of a box in each axis be given by ${x_k}^{-}, {x_k}^{+}$, where $k = 1,...,M$ and $x \in \reals^M$. Here, ${x_k}^{-} \le x_k \le {x_k}^{+}$. Thus, an aligned box is denoted by $B({x}^{-},{x}^{+})$ (or, succinctly by $B$). 

To ensure we obtain homothets, define a prototype box, $\hat{B}^{0}$. The choice of $\hat{B}^{0}$ can be arbitrary; for example, one could consider a square in $\reals^2$, or a hypercube in $\reals^M$. In our case, we find $\hat{B}^{0}$ by solving for the largest volume box \cite{VaBo1998} that fits in a representative polytope, $\cP^0$, chosen from the $n_d$ given polytopes. Let the distances of the edges of $\hat{B}^{0}$ be, $d_k^{0} = ({x}_k^{+}-{x}_k^{-}), \; k = 1,2,...,M$. Then, the ratios of distances w.r.t $d_1^{0}$ are, ${r}_{1,k}^{0} = \frac{d_1^{0}}{d_k^{0}}. \; k = 2,..., M$. From here on, we require all boxes, $B$, must be homothets of $\hat{B}^0$. 

Given $\cP$, in order to find a homothet, $B$, with maximum volume, and ensure $B \subseteq \cP$, the following problem must be solved, 
\begin{subequations}
\begin{align}
&(P2) \;  \max_{{x}^{+},{x}^{-}} \;  \prod_{k=1}^{M} \big( {x}_k^{+} - {x}_k^{-}  \big)    \label{eq:LB1}\\
& \quad \text{s.t.}  \; A^{+}{x}^{+} -  A^{-} {x}^{-} \le b ,   \label{eq:LB2} \\
& \quad {x}_k^{-} \le {x}^{+}_k, \quad k = 1,2,...,M, \label{eq:LB3} \\
& \quad \big({x_1}^{+} - {x_1}^{-}\big) = {r}_{1,k}^{0} \big({x_k}^{+} - {x_k}^{-}\big), \; k = 2,...,M, \label{eq:LB6}
\end{align}
\end{subequations}
where $A^{+}_{ij} = \textrm{max}\{0,A_{ij} \}$ and  $A^{-}_{ij} = \textrm{max}\{0,-A_{ij} \}$, with $i,j$ being the row and column indices of $A^{+}, A^{-}$ \cite{VaBo1998}. Note that the objective function \eqref{eq:LB1} can be replaced by $\sum_{k=1}^{M} \textrm{log} \big( {x}^{+} - {x}^{-}  \big)$, which will be a convex problem \cite{Muller2017}. Constraint \eqref{eq:LB6} ensures $B$ will be a homothet of $\hat{B}^0$. Next, we show how (P2) can be used in a multi-stage algorithm for decomposing $\cP$ into a number of homothets. 

\begin{figure}[t!]
	\begin{center}
		\includegraphics[scale=0.31]{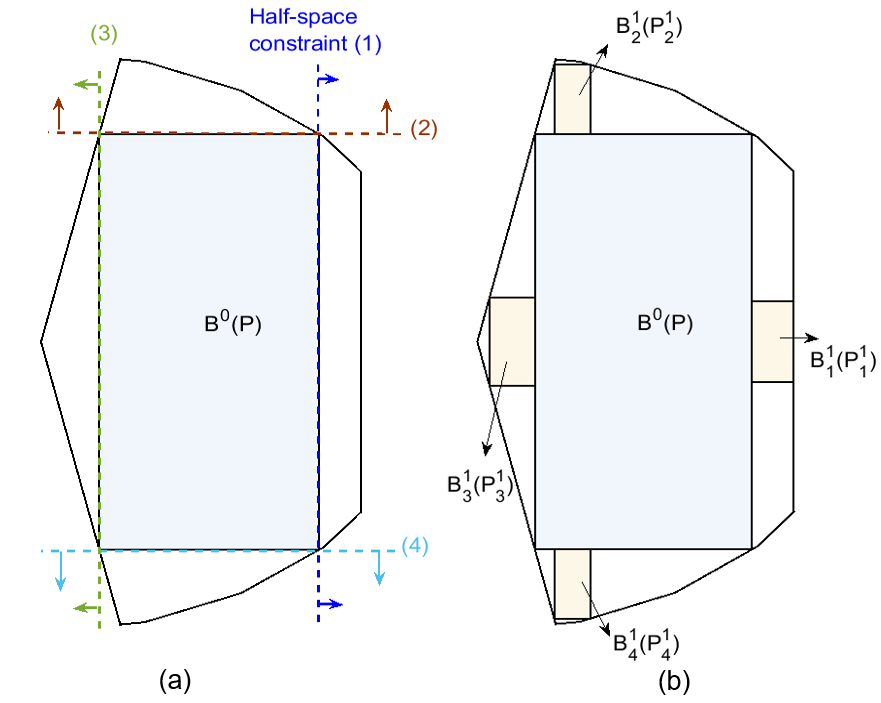}
		\caption{(a) Homothet-based polytope decomposition: $B^0(\cP)$ obtained in $s=0$ and the half-space constraints originating from $B^0(\cP)$. (b) Decomposition after completion of stages, $s=0$ and $s=1$. }
		\label{fig:ConceptSPF0}
	\end{center}	
\end{figure}

Let $n_s$ represent the total number of stages and $s = 0, 1, ..., n_s$ denote the stage index. At $s=0$, (P2) is solved for $\cP$ to obtain $B^{0}(\cP)$, i.e. $B$ is an outcome of the polytope. Here, $B^{0}(\cP) = \beta_{0} \hat{B}^{0} + c_{0}$, i.e. a homothet of $\hat{B}^{0}$. Fig. \ref{fig:ConceptSPF0}(a) shows $B^{0}(\cP)$ inside $\cP$. 

Next, at $s=1$, additional homothets will be sought in each region outside $B^{0}(\cP)$, but inside $\cP$. In general, since each $B$ in $\reals^M$ has 2M half-space constraints, let $\sigma = 1, 2,...,2M$ be their index. As shown in Fig. \ref{fig:ConceptSPF0}, $B^{0}(\cP)$ is defined by four half-space inequalities in $\reals^2$. Each region outside $B^{0}(\cP)$, but inside $\cP$, can be defined using the half-space inequalities of $B^{0}(\cP)$, except the sign of the inequalities must be reversed, as illustrated in Fig. \ref{fig:ConceptSPF0}(a). By construction, each region outside $B^0(\cP)$, but inside $\cP$ is convex and compact. Let $\cP^{s}_{\sigma}$ denote the updated polytope corresponding to half-space inequality $\sigma$ and $B^{s}_{\sigma}(\cP^{s}_{\sigma})$ be the corresponding solution obtained by solving (P2). 
The HPD concept is further illustrated in Fig. \ref{fig:ConceptSPF0}(b) with stages $s=0$ and $s=1$ solved. The decomposition can continue up to $n_s$ stages. The algorithm is summarized below.


  %
\vspace{2pt}

\noindent
\textbf{Algorithm (HPD):}

\noindent
\textbf{S1.} $s=0$. Starting with $\cP$, compute $B^0(\cP)$. \\
\textbf{S2.}
$s = s+1$. If $s \le n_s$ proceed, else stop.\\
\indent	\textbf{S2a.}  $s=1$, $\sigma = 1$. Define $\cP_1^1 = \cP \cap$ half-space inequality $\sigma=1$. Compute $B_1^1 (\cP_1^1)$.\\
\indent	\textbf{S2b.} $s=1$, $\sigma = 2$. Define $\cP_2^1 = \cP \cap$ half-space inequality $\sigma=2$. Compute $B_2^1 (\cP_2^1 )$.\\
\vdots

(Continue till $s=1$, $\sigma = 2M$.)\\
\textbf{S3.}  $s = s+1$. If $s \le n_s$, proceed else stop.\\
\indent	\textbf{S3a.} $s = 2$, $\sigma = 1$. Define the new polytope $P^2_{1,1} = \cP_1^1 \cap$ half-space inequality $\sigma=1$. Compute $B^2_{1,1} ( \cP^2_{1,1})$.\\
\vdots

(Continue till $s=2$, $\sigma = 2M$.)\\
\vdots

\noindent
\textbf{S4.} Continue until $s > n_{s}$.

\subsection{Convergence of Polytope Decomposition and M-Sum}\label{sec:HPDConvLemmas}

Following the HPD algorithm, at every stage, new regions will be covered by solving (P2), unless the entire polytope has already been covered. 

\begin{observation}\label{LemmaHPD1}
If at the end of stage $s$, $\big(\cP - \cup_{\forall s,\sigma} B^s_{\sigma}(\cP^s_{\sigma})\big) \ne \emptyset$, then $\textrm{vol}\big(B^{s+1}_{\sigma}(\cP^{s+1}_{\sigma})\big) > 0$.	
\end{observation}

Let $B^{s}_{\tilde{\sigma}}(P^{s}_{\tilde{\sigma}})$ be any box obtained at stage $s$ that follows from constraint $\tilde{\sigma}$ of $B^{s-1}_{{\sigma}}(P^{s-1}_{{\sigma}})$. Then, the following holds.

\begin{observation}\label{LemmaHPD2}
$\textrm{vol}\big( B^{s}_{\tilde{\sigma}}(\cP^{s}_{\tilde{\sigma}}) \big) \le \textrm{vol}\big(B^{s-1}_{\sigma}(\cP^{s-1}_{\sigma})\big)$, $\forall \tilde{\sigma}$.
\end{observation}

The above holds because otherwise it would contradict the solution of (P2). After completion of the decomposition phase for $n_d$ polytopes, the approximate M-sum polytope, $\cP^{\textrm{Aggr}}$, can be obtained using \eqref{eq:MsumUnionHomothet}. Then, the following asymptotic result holds.

\begin{proposition}\label{ThmUnionMsum}
	$\cP^{\textrm{Aggr}} \rightarrow \cX^{\textrm{Aggr}}$ as $s \rightarrow \infty$.
\end{proposition}

\begin{remark}
The HPD algorithm and the union-based M-sum computation procedure using the axis-aligned boxes are general to any dimension. 
\end{remark}




\subsection{Practical Considerations}\label{sec:Practical}

While the proposed algorithm can guarantee asymptotic convergence to the true M-sum, considering a large number of sub-polytopes for $n_d$ devices can be computationally challenging. Hence, a number of strategies can be considered. 

The HPD algorithm, presented in section \ref{sec:PDecomposition}, stops after completing $n_{s}$ stages. Alternatively, the stopping condition can be based on the volume of $B^{s}_{\sigma}$. After reaching a certain threshold, one can stop because all subsequent boxes will be smaller (by Observation \ref{LemmaHPD1}). 

As discussed earlier, for many applications, it could be sufficient to utilize \eqref{eq:optUnion}. Hence, instead of computing the entire M-sum, one can use a set of candidate polytopes and solve \eqref{eq:optUnion}. Section \ref{sec:Results} provides detailed examples on how to efficiently choose such candidates. In $\reals^2$, one can also consider computing the convex hull (C-hull) of the aggregate boxes to obtain a single M-sum approximation polytope, which will also be shown in section \ref{sec:ResultsINV}. 

Section \ref{sec:Results} provides numerical examples, considering inverter polytopes in $\reals^2$ and storage-like loads in $\reals^6$, to analyze the performance of the proposed schemes and discusses trade-offs. 

\subsection{Simplification for Axis-Aligned Boxes}\label{sec:Interval}

The choice of axis-aligned boxes leads to a further simplification due to applicability of interval analysis techniques \cite{Moore2009}. For the $i$-th DER, define an interval $I({x}_{i,k}^{-},{x}_{i,k}^{+})$, $k = 1,2,...M$, such that ${x}_{i,k} \in I({x}_{i,k}^{-},{x}_{i,k}^{+})$ implies ${x}_{i,k}^{-} \le {x}_{i,k} \le {x}_{i,k}^{+}$. Then, $\sum_{i=1}^{n_d} {x}_{i,k}^{-} \le {x}^{\textrm{Aggr}}_{k} \le \sum_{i=1}^{n_d}{x}_{i,k}^{+}$, where ${x}^{\textrm{Aggr}}_{k}=\sum_{i=1}^{n_d} {x}_{i,k}$ for all ${x}_{i,k} \in I({x}_{i,k}^{-},{x}_{i,k}^{+})$, $i = 1,2,...,n_d$ and $k = 1,2,...M$. In $\reals^M$, we obtain $B(\sum_{i=1}^{n_d} {x}_{i,k}^{-},\sum_{i=1}^{n_d} {x}_{i,k}^{+})$, which is by default an inner approximation to $\cX^{\textrm{Aggr}}$. The decomposition procedure described in HPD remains exactly the same, except, we can relax \eqref{eq:LB6}. The convergence results discussed in section \ref{sec:HPDConvLemmas} trivially extend to the case of applying interval arithmetic on axis-aligned boxes. 


\section{Numerical Results}\label{sec:Results}

\subsection{Performance of Union-based M-sum for Inverters}\label{sec:ResultsINV}

Consider four different inverters, $\cX_i(S_i,\overline{P_i},\theta_i), i = 1, 2, 3, 4$,  with parameters (A) (1, 0.9,$\pi/2$ rad), (B) (1, 0.8, 1.37 rad), (C) (1, 0.6, 1.37 rad), and (D) (1, 0.3,$\pi/2$ rad).

\begin{figure}[t!]
	\begin{center}
		\includegraphics[scale=0.35]{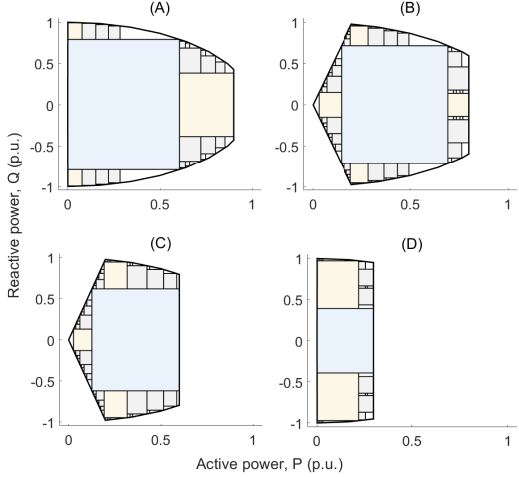}
		\caption{Decomposition of four inverter polytopes using the homothet-based polytope decomposition algorithm.}
		\label{fig:SPF4inv}
	\end{center}
    \vspace{-.4cm}
\end{figure}

The results obtained by applying the HPD algorithm to each of the four inverters are shown in Fig.~\ref{fig:SPF4inv}. The HPD was solved using CVX \cite{CVX}. All computations were performed on a computer with Intel Core i-5 3.20 GHz processor with 8 GB RAM.

\begin{figure}[b!]
	\begin{center}
		\includegraphics[scale=0.35]{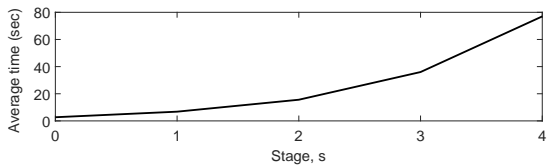}
		\caption{Average time for decomposition up to stage $s$.}
		\label{fig:time1}
	\end{center}
    \vspace{-.4cm}
\end{figure}

The area ratios (approximated area divided by the area of the $i$-th inverter polytope) after completion of each stage are given in Table~\ref{tab:performanceDecompose}. From the area ratios, on average, the total area captured after stage, $s = 0$ was 58 \%, whereas after $s = 4$ was 95 \%. The average time taken to complete decomposition up to a each stage $s$ is shown in Fig. \ref{fig:time1}.  Up to $s = 1$, the average computation time was only 6.7s, whereas the area ratios averaged at 79\%, a 21 \% increase from the case of $s = 0$. 


\begin{table}[b!]
	\caption{Area covered, as a fraction of the area of true M-sum, after every stage of HPD, for each inverter.} \label{tab:performanceDecompose}
	\centering
	\begin{tabular}{|c|c|c|c|l|} \hline
		Stage, $s$ & (A) 	& (B) & (C) & (D) 	\\ \hline
		 0 	  & 0.64 	& 0.65 	& 0.64  & 0.40	\\
		 1 	  & 0.81 	& 0.74 	& 0.76 	& 0.84	\\
		 2    & 0.87 	& 0.85 	& 0.86  & 0.89	\\ 
		 3    & 0.92 	& 0.91 	& 0.92 	& 0.95	\\ 	
		 4    & 0.94 	& 0.93 	& 0.96 	& 0.96	\\ 	\hline
	\end{tabular}
\end{table}

\begin{figure}[t!]
	\begin{center}
		\includegraphics[scale=0.33]{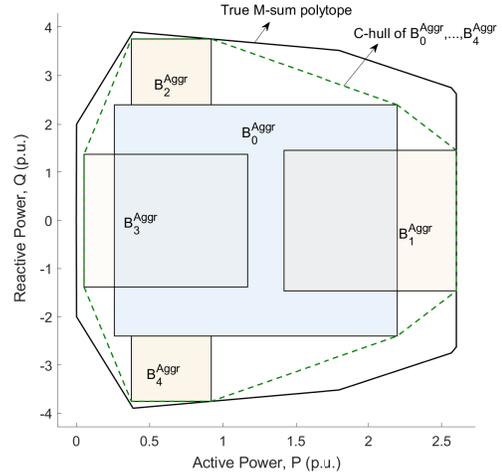}
		\caption{M-sum approximation using $B^{Aggr}_0$ and $B^{Aggr}_{\sigma}$, $\sigma=1,...,4$, the C-hull of the aggregate boxes, and the true M-sum polytope,  $\cP^{\textrm{Aggr}}$. }
		\label{fig:HeteroSumCase1}
	\end{center}	
\end{figure}

The result of the aggregation is shown in Fig.~\ref{fig:HeteroSumCase1}. As discussed in sections \ref{sec:DistributivityMsum} and \ref{sec:Practical}, instead of computing the entire M-sum, we consider a limited set of \textit{candidate} boxes, which in this case correspond fist selecting only stage 0 and stage 1 boxes for each polytope. Secondly, for computing the M-sums, instead of considering all combinations of unions, consider : $B^{Aggr}_0 = \oplus_{i=1}^{n_d}B^0_i(\cP_i)$, 
$B^{Aggr}_{\sigma} = \oplus_{i=1}^{n_d} B^1_{i,\sigma}(\cP^{1}_{\sigma,i}), \; \sigma=1,2,3,4$. In case, any $B^1_{i,\sigma}(\cP^1_{\sigma,i})$ is degenerate, it was replaced with $B^0_{i}(\cP_i)$. These five aggregate boxes are shown in Fig.~\ref{fig:HeteroSumCase1}. 

The actual M-sum polytope, also shown Fig.~\ref{fig:HeteroSumCase1}, was obtained using the MPT toolbox \cite{MPTtoolbox}. Finally, the C-hull of these boxes was computed in MATLAB and is shown in Fig. \ref{fig:HeteroSumCase1}. Since all vertices of the aggregate boxes lie inside the true M-sum, which is convex and compact, the C-hull of these aggregate boxes is also an inner approximation to the true M-sum polytope. 

To assess the M-sum approximation accuracy, we computed the area ratios for our approximated boxes and the C-hull and compared these against the area of the true M-sum polytope. Using only $B^{Aggr}_0$, the M-sum approximation accuracy was 52\%. Using both stage 0 and stage 1 candidate boxes, the accuracy increased to 71\%, thus demonstrating the effectiveness of considering multiple homothets per device. Finally, with C-hull, the accuracy was 85\%.

As discussed before, considering all combinations of unions would cause exponential growth in complexity \cite{Weibel2007}. Instead, our policy used a fixed set of candidates. Once the results of the decomposition were available, computing the five aggregate boxes was a trivial task. Furthermore, given that the vertices of the aggregate boxes were readily available, the C-hull was also inexpensive to compute in $\reals^2$, where C-hull computational has complexity O($n_v$ log $n_v$) with $n_v$ being the number of vertices \cite{AvisHULL}. .

\subsection{Performance of Union-based M-sum for Controllable Loads and Storage Devices}

To evaluate the performance of our proposed scheme in $\reals^M$, consider storage-like loads and a 3 hour planning horizon with 30 minutes time steps. Hence, $M$ = 6. Consider $n_d$ = 100 devices with parameters, $\underline{P} = 0$, $\overline{P} \in [3,4.5]$ kW,  
$a \in [0.9,1]$, $e_{0} \in [0.2,0.6]$ (normalized), and $\gamma \in [0.035,0.053]$. 
First, the polytope decomposition algorithm of \ref{sec:PDecomposition} was applied, with and without \eqref{eq:LB6} active. Imposing \eqref{eq:LB6} generally led to slower coverage of the feasible region of $\cP$. Hence, we relaxed \eqref{eq:LB6} in $\reals^6$ to improve volume coverage per stage. To evaluate accuracy of the decomposition phase, volume ratios were considered. Since volume computation in high dimension is an NP-hard problem, this is done here by finding a bounding box and generating Monte Carlo samples \cite{Barot2017}. The volume ratios averaged 56\% at the end of $s=0$, 67\% at the end of $s = 1$ and 82\% at the end of $s = 2$. The decomposition up to $s = 1$, for each polytope, took on average 6.2s. Recall that the decompositions can be performed in parallel.

For computing the M-sum, we consider a limited set of candidate boxes from $s = 0$ and $s = 1$. Following the same procedure in section \ref{sec:ResultsINV}, we obtain $(2M+1) = 13$ boxes in $\reals^6$. 
In order to compute true M-sum polytope, we again used the MPT toolbox \cite{MPTtoolbox}. However, using MPT, it was only possible to obtain M-sums of 5 devices at a time, beyond which no solutions were reached in several hours. Hence, we sampled 5 devices from the population randomly and repeated the volume computation to obtain an estimate of the volume ratios. The accuracy of the M-sum approximation, compared to the actual volume, was 44\% using only stage 0 boxes, and increased to 74\% using the 13 candidate boxes. In our case, since the M-sum computation simply required an algebraic summation of the upper and lower bounds of intervals (or the scaling and translation coefficients of homothets in the previous example), the M-sum computation time was negligible and did not depend on the population size. 


\begin{table}[t!]
	\caption{Numbers of Convex Problems (P2) solved by end of stages, $s = 0, 1, 2, 3$, with increase in dimension, $M$} \label{tab:cvxDecompose}
	\centering
\begin{tabular}{|c|ccccccc|}
\hline
 &  &  &   & $M$ &  & &     \\
Stage, $s$ & 2 & 3 & 4 & 5 & 6 & 7 & 8  \\
 \hline
 0 & 1 & 1 & 1 & 1 & 1 & 1 & 1  \\
 1 & 5 & 7 & 9 & 11 & 13 & 15 & 17   \\
 2 & 21 & 43 & 73 & 111 & 157 & 211 & 273 \\
 3 & 85 & 259 & 585 & 1111 & 1885 & 2955 & 4369  \\
\hline
\end{tabular}
\end{table}

Given that at each decomposition stage (P2) must be solved, we can express time complexity in terms of the number of convex problems solved. During decomposition with axis-aligned boxes, two half-space inequalities are considered in each dimension. Given $x \in \reals^M$, each box is defined by $2M$ constraints. Then, at stage $s$, for each box, $2M$ additional constraints are introduced for the subsequent stage. Thus, at stage, $s$, (P2) is solved at most $(2M)^s$ times. Table \ref{tab:cvxDecompose} shows the maximum number of times (P2) must be solved by the end of stage $s$. While choosing a higher $s$ assists in achieving convergence to the true M-sum, the complexity grows exponentially. Hence, in our approach, in $\reals^4$ and above we suggest to choose $s = 1$, which ensures a polynomial time algorithm. Using our of candidate selection for computing a finite set of aggregate boxes, in $\reals^M$, we obtain exact $(2M+1)$ boxes. For the purpose of solving optimization problems, recall that (P1) can now be solved subject to each of these $(2M+1)$ boxes, in parallel. Since only $2M$ half-space constraints are required to represent the boxes, the optimization problem is generally significantly simpler than (P0) subject to the true M-sum polytope.

\section{Conclusions}\label{sec:Conclusion}
In this paper, we developed and compared algorithms to compute inner approximations of the Minkowski sum of convex polytopes. As an application, we considered the computation of the feasibility set of aggregations of distributed energy resources (DERs), such as solar photovoltaic inverters, controllable loads, and storage devices. A convex polytopic representation for a feasible operating region of inverter interfaced DERs was developed first. We showed how homothets can be used to compute the M-sum and obtained analytical expressions in special cases. However, as heterogeneity increases, using a single homothet per device, may result in highly conservative inner approximation of the M-sum. Hence, to fully account for the heterogeneity in the DERs while ensuring an acceptable approximation accuracy, we leveraged  a union-based computation that advocates a homothet-based polytope decomposition. We show that the proposed algorithm can guarantee the inner approximation asymptotically converges to the true M-sum. However, union-based approached can in general lead to high-dimensionality concerns; to alleviate this issue, this paper shows how to define candidate sets to reduce the computational complexity. Accuracy and trade-offs have been analyzed through numerical examples. The flexibility polytopes of inverter-interfaced devices, controllable loads and storage can be integrated in power systems planning tools to provide various power system services. 



%

\appendix

\subsection{Proof of Proposition 1}

Consider $\cX(S,\underline{P},\overline{P})$. Assume $\underline{P} = -1$ and $\overline{P} = 1$. 
Then, the area ratio, $\eta$, can be expressed as, 
\begin{equation}\label{eq:polyratio}
\eta = \frac{A_{\cP}}{A_{\cX}} = \frac{0.5 N S^2 \text{sin}(2\pi/N)}{\pi S^2} =\frac{sin(\alpha)}{\alpha}. 
\end{equation} 
where $\alpha = 2\pi/N$ and since we chose $N \ge 4$, $0 \le \alpha \le \pi/2$. 
As $N$ increases, $\alpha$ decreases. Then, by application of the L'Hopital's rule, $\underset{\alpha \rightarrow 0}\lim \frac{sin(\alpha)}{\alpha}= \underset{\alpha \rightarrow 0}\lim cos(\alpha) = 1$. Therefore, $\eta \rightarrow 1$ as $\alpha \rightarrow 0$ (i.e. $N \rightarrow \infty$).


When considering $\cX(S,\overline{P},\theta)$, with $\overline{P} = 1,\theta = \pi/2$, the third fraction in \eqref{eq:polyratio} will additionally have 0.5 multiplied both in numerator and denominator (due to half-circle), which cancel out. Thus, the ratio $\eta = \frac{sin(\alpha)}{\alpha}$ holds, and again $\eta \rightarrow 1$ as $\alpha \rightarrow 0$.$\qed$

\subsection{Proof of Proposition 2}
\begin{figure}[b!]
	\begin{center}
		\includegraphics[scale=0.23]{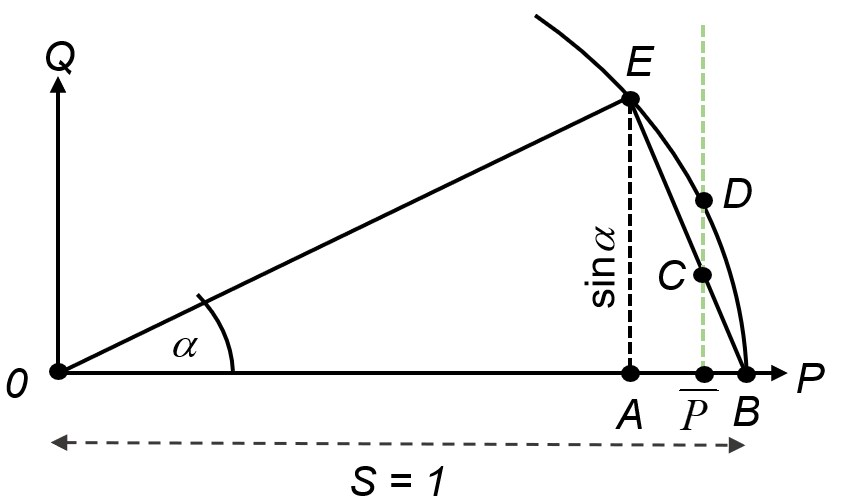}
		\caption{Discretization of circle and adding constraint.}
		\label{fig:prop1}
	\end{center}	
\vspace{-.4cm}
\end{figure}

Without loss of generality, assume $S=1$, as shown in Fig. \ref{fig:prop1}. $\alpha$ is the angle formed by fitting an N-sided polygon inside $P^2+Q^2 = 1$. Thus, the area of BCED $=$ (area of sector OBDE - area of $\Delta$OBE), is the approximation error. 

The area of OBDE $= \frac{\alpha}{2}$ and the area of $\Delta$OBE $= \frac{\sin\alpha}{2}$. By inclusion, $\frac{\alpha}{2} \ge \frac{\sin\alpha}{2}.$

However, for $\alpha$ sufficiently small, $\alpha \approx \sin \alpha$ (small angle approximation). Thus, $\frac{\alpha}{2} \approx \frac{\sin\alpha}{2}$, i.e. area of BCED $\approx$ 0.

Consider $\underline{P} = -1$ and $\overline{P} < 1$. Assume $\overline{P}$ lies between A and B. Now, the area of CED is the approximation error when the sector OBE is approximated by $\Delta$OBE and constraint $P \le \overline{P}$ is added. For $\overline{P} < 1$, area of CED $<$ the area of BCED. Hence, with $\alpha$ sufficiently small, area of BCED $\approx$ 0, hence, the area of CED $\approx$ 0. We can show the same considering $\underline{P} > -1$. Combining these, we obtain, $A_{\cP} \approx A_{\cX}$, i.e. $\eta \approx 1$, for $N$ sufficiently large.

It is also easy to verify that the same holds when considering $\cX(S,\overline{P},\theta)$, with $\overline{P}<1,0 \le \theta < \pi/2$. $\qed$

 
\subsection{Proof of Theorem 1}

For two inverters, assume $S_1 \neq S_2$, $\underline{P_1} = \underline{P_2} = \underline{P_0}$ and $\overline{P_1} = \overline{P_2} = \overline{P_0}$. Let $m_{2,1} > 0$ be the ratio of the rated powers of the two inverters, i.e. $m_{2,1} = S_2/S_1$. Choose prototype set, $\cX_{0} =  {\cX_1}$. Then, ${\cX_2}= m_{2,1} \cX_0$. By \eqref{eq:M-sumhomothet2}, we can write, 

$$ ({\cX_1} \oplus {\cX_2}) = (1+ m_{2,1})\cX_{0}.$$
Similarly for $n_d$ devices, with $\cX_{0} = \cX_{1}$ we obtain, 
$${\cX}^{\textrm{Aggr}} = (1+ m_{2,1}+...+m_{n_d,1})\cX_{0}.$$
where $m_{i,1} = S_i/S_1, i = 2, ..., n_d$.

Thus, plugging in $m_{i,1}$, we verify that  $\cX \big(\sum_{i=1}^{n_d} S_i, \underline{P_0},\overline{P_0} \big) = {\cX}^{\textrm{Aggr}}$.$\qed$

\subsection{Proof of Theorem 2}

Consider $n_d$ devices with feasible sets $\cX_i \big(S_i, \overline{P_i}, \theta_i \big)$. Let $S_0 = \min_i \big(S_i), \overline{P_0} =  \min_i \big(\overline{P_i}, \big), \theta_0 = \min_i \big(\theta_i \big)$ and obtain $\cX_0(S_0, \overline{P_i}, \theta_0)$. 
In this case, we obtain $\cX_0 = \cX_1 \cap \cX_2  \cap ...  \cap \cX_{n_d}$. Hence, $\cX_0 \subseteq \cX_i, \forall i$.
Approximate every $\cX_i$ by $\cX_0$, which gives  $\beta_i = 1$, $t_i = 0$ in \eqref{eq:M-sumhomothet2}. Hence, $\sum_{i=1}^{n_d} \beta_i = n_d$. Therefore,
\begin{equation}\nonumber
\cX \big(n_d S_0, \overline{P_0}, \theta_0 \big) = n_d \cX_0 \subseteq {\cX}^{\textrm{Aggr}}. 
\end{equation}
In the homogeneous case, $\cX_0 = \cX_i, \forall i$, thus strict equality holds.
$\qed$

\subsection{Proof of Proposition 3}

To prove, first we need to show, as $s \rightarrow \infty$, $\cP_i \rightarrow \cX_i$. 
Assume $\textrm{vol}(\cP_i) = \textrm{vol}\big(\cup_{\forall s,\sigma} B^{s}_{\sigma}(\cP^{s}_{\sigma}) \big) + \delta$, where $\delta \in \reals$. First consider $\delta > 0 $. This implies $\big(\cP - \cup_{\forall s,\sigma} B^s_{\sigma}(\cP^{s}_{\sigma})\big) \ne \emptyset$, and for at least one $B^{s}_{\sigma}(\cP^{s}_{\sigma})$, we can find a region $\cP^{s+1}_{\tilde{\sigma}}$, where applying HPD results in $\textrm{vol}\big(B^{s+1}_{\tilde{\sigma}}(\cP^{s+1}_{\tilde{\sigma}})\big) > 0$. But this is a contradiction since Observation \ref{LemmaHPD1} and \ref{LemmaHPD2} imply as $s \rightarrow \infty$, $\textrm{vol}\big(B^{s}_{\sigma}(\cP^{s}_{\sigma})\big) \rightarrow 0, \forall s, \forall \sigma)$. Next consider $\delta < 0$. Then, $\textrm{vol}\textrm{vol}\big(\cup_{\forall s,\sigma} B^{s}_{\sigma}(\cP^{s}_{\sigma}) \big) > \textrm{vol}(\cP_i)$, but this is again a contradiction since $\cup_{\forall s,\sigma} B^{s}_{\sigma}(\cP^{s}_{\sigma}) \subseteq \cP_i$. Hence, $\delta$ must be equal to 0. Hence, $s \rightarrow \infty$, $\cP_i \rightarrow \cX_i$.

Similarly, from \eqref{eq:MsumUnionHomothet}, it follows that, $s \rightarrow \infty$, $\cP^{\textrm{Aggr}} \rightarrow \cX^{\textrm{Aggr}}$. $\qed$



\bibliographystyle{IEEEtran}
\bibliography{biblio.bib}


\end{document}